\documentclass{ws-mpla}
\usepackage{color}
\usepackage{hyperref}
\usepackage{mathtools}
\usepackage{url}
\hypersetup{colorlinks,urlcolor=cyan,citecolor=green,linkcolor=blue,filecolor=black}
\begin{document}

\catchline{}{}{}{}{}
\title{Majorana transformation of the Thomas–Fermi equation demystified}

\author{Abdaljalel Alizzi}
\address{Novosibirsk State University, 630 090, Novosibirsk, Russia. \\ abdaljalel90@gmail.com}
\author{Zurab K. Silagadze}
\address{Budker Institute of Nuclear Physics and \\ Novosibirsk State University, 630 090, Novosibirsk, Russia. \\ silagadze@inp.nsk.su}

\maketitle

\pub{Received (Day Month Year)}{Revised (Day Month Year)}

\begin{abstract}
The Majorana transformation makes it possible to reduce the Thomas-Fermi equation to a first-order differential equation. This reduction is possible due to the special scaling property of the Thomas-Fermi equation under homology transformations. Such reductions are well known in the context of stellar astrophysics, where the use of homology-invariant variables has long proved useful. We use homology-invariant variables in the context of the Thomas-Fermi equation to demystify the origin of the otherwise mysterious Majorana transformation.

\end{abstract}

\section{Introduction}
Llewellyn H. Thomas and Enrico Fermi independently proposed the semi-classical statistical model of multi-electron atoms in 1926  \cite{thomas_1927,Fermi_1927}. The fundamental concept of the model is to consider the electron cloud around the nucleus as a degenerate Fermi-Dirac fluid,  which is held in a hydrostatic equilibrium at absolute zero temperature in a self-consistent electric field. The Hartree-Fock self-consistent-field model can be approximated by the Thomas-Fermi method in a semi-classical manner, allowing for the introduction of numerous improvements of the Thomas-Fermi model \cite{Kompaneets_1956}. 

The local electron density is the key variable in the Thomas-Fermi technique, which completely omits the consideration of the multi-electron wave function. In this regard, the Thomas-Fermi model is the simplest density-functional theory \cite{Morgan_2006}. As a result of the development of the Thomas-Fermi model, culminating in the seminal papers of Hohenberg and Kohn \cite{Hohenberg_Kohn} and Kohn and Scham \cite{Kohn_Sham}, the modern density-functional theory has emerged, which has become a highly effective tool in chemistry and materials science \cite{Burke_2012,Jones_1989}.

It is not unexpected that the elegant and straightforward Thomas-Fermi model has inspired a substantial body of literature. We merely draw attention to a few reviews \cite{March_1957,Spruch,Parr_Weitao} and textbook expositions \cite{Bethe_Jackiw,Flugge,Kleinert,Landau_1991,Budker_2008}. 

Majorana discovered an intriguing transformation of the Thomas-Fermi equation in 1928, but it was never published and was overlooked until recently \cite{Esposito_2002,DiGrezia_2004,Esposito_2015}. The Majorana transformation, whose origin appears to be rather enigmatic, enables one to reduce the degree of Thomas-Fermi equation and this is useful in the geometric theory of this equation \cite{Seiler_2023}. Further study shows that this reduction is possible because the Thomas-Fermi equation is scaling-invariant, and Majorana's method can be applied to other ordinary differential equations that satisfy Majorana-type scaling relations \cite{Esposito2002}. The Majorana transformation appears less perplexing also from the perspective of Lampariello dynamical systems \cite{Rosu_2017}. But like the fate of its renowned author \cite{Esposito_2017}, there is still a certain air of mystery surrounding this transformation. 

To our amazement, or perhaps there is no surprise here, as we are not the first to find inspiration \cite{Eggleton_2011} in the masterpiece of Chandrasekhar \cite{Chandrasekhar}, the insight came from a completely unanticipated source, namely from the study of the structure of stellar bodies.

Almost immediately after the appearance of the Thomas-Fermi equation, Milne noticed that it was similar to the Lane-Emden equation used to study the polytropic equilibrium of a star \cite{milne_1927}. The reduction of this equation to a single first-order equation was already obtained by Emden, and the possibility of such a reduction is closely related to the invariance of the Lane-Emden equation under homology transformations \cite{Hopf,Chandrasekhar}.

In this short note, we argue that the invariance of the Thomas-Fermi equation under homology transformations is the key to the riddle of the Majorana transformation, allowing it to be fully elucidated.

\section{Homology-invariant variables and the reduction of the Emden-Fowler equation}
Let us explain the main idea using the example of the Emden-Fowler equation in the form (this particular case of the Emden-Fowler equation is essentially equivalent to the Lane-Emden equation) \cite{Hille1970}
\begin{equation}
    y^{\prime\prime}(x)=x^{1-p}y^p(x),
    \label{eq1}
\end{equation}
where prime denotes differentiation with respect to $x$, and $p$ is some real number. The general solution of this equation, which is a second-order differential equation, must be characterized by two integration constants. However, the equation (\ref{eq1}) remains invariant under the homology transformations
\begin{equation}
    x\to \lambda x,\;\;y\to\lambda^{-q}y,\;\;q=\frac{3-p}{p-1}.
    \label{eq2}
\end{equation}
In other words, if $y(x)$ is a solution of (\ref{eq1}), then for any real number $\lambda$ the homologically related function $y_\lambda(x)=\lambda^{-q }y(\lambda x) $ is also a solution of the same equation. Then one constant of integration must be trivial in the sense that it simply determines the scaling factor $\lambda$, and it must be possible to bring the equation to the first order \cite{Chandrasekhar,Horedt}.

In fact, there are infinitely many ways in this case to reduce to a first-order differential equation, since all we have to do is choose homology-invariant variables \cite{Chandrasekhar,Horedt}. For example, we can use Coppel's homology-invariant variables \cite{Hille1970,Coppel}
\begin{equation}
    v=\frac{xy^\prime}{y},\;\;\;u=\frac{x^{2-p}y^p}{y^\prime},\;\;\;uv=x^{3-p}y^{p-1}.
    \label{eq3}
\end{equation}
Evidently, $x\frac{du}{dx}$ and $x\frac{dv}{dx}$ are by themselves homology invariant quantities and we expect them to be expressible in terms of $v$ and $u$. Indeed, this is confirmed by simple calculations (in which we make use of (\ref{eq1})):
\begin{equation}
    x\frac{dv}{dx}=x\left [\frac{y^\prime}{y}+\frac{xy^{\prime\prime}}{y}-\frac{xy^{\prime\,2}}{y^2}\right ]=x\frac{y^\prime}{y}+x^{3-p}y^{p-1}-x^2\frac{y^{\prime\,2}}{y^2}=v+uv-v^2,
    \label{eq4}
\end{equation}
and
\begin{eqnarray} &&
   x\frac{du}{dx}=x\left [(2-p)\frac{x^{1-p}y^p}{y^\prime}+px^{2-p}y^{p-1}-\frac{x^{2-p}y^p}{y^{\prime\,2}}y^{\prime\prime}\right ]= \nonumber \\ &&
   (2-p)\frac{x^{2-p}y^p}{y\prime}+px^{3-p}y^{p-1}-\frac{x^{2(2-p)}y^{2p}}{y^{\prime\,2}}=(2-p)u+puv-u^2.
    \label{eq5}  
\end{eqnarray}
Therefore, the first order equation equivalent to the Emden-Fowler equation (\ref{eq1}) is
\begin{equation}
    \frac{du}{dv}=\frac{u\left(2-p+pv-u\right)}{v\left(1+u-v\right)}.
\label{eq6}
\end{equation}
Any two independent homology-invariant variables will do a similar job of reducing the Emden-Fowler equation to an equivalent first order equation. For example, we can use the Milne homology-invariant variables, which are well known in astrophysics (keep in mind,  to avoid having a multitude of variables, we use the same notations $u$ and $v$, although now they have a different meaning) \cite{Chandrasekhar,Horedt}:
\begin{equation}
    v=\frac{x\theta^\prime}{\theta},\;\;\;u=\frac{x\theta^p}{\theta^\prime},\;\;\;uv=x^2\theta^{p-1},
    \label{eq7}
\end{equation}
where $\theta=y/x$, and in terms of $\theta$ equation (\ref{eq1}) becomes essentially Lane-Emden equation (in stellar astrophysics, the $\theta^p$ term in the right-hand-side of the Lane-Emden equation has the opposite sign)
\begin{equation}
    \theta^{\prime\prime}+\frac{2}{x}\theta^\prime=\theta^p.
    \label{eq8}
\end{equation}
A simple calculation will give
\begin{equation}
    x\frac{dv}{dx}=v(u-v-1),\;\;\;x\frac{du}{dx}=u(3+pv-u),
    \label{eq9}
\end{equation}
and so in this case the first order equation has the form
\begin{equation}
    \frac{du}{dv}=\frac{u(3+pv-u)}{v(u-v-1)}.
    \label{eq10}
\end{equation}

\section{Majorana homology-invariant variables}
The Thomas-Fermi equation is a special case of the Emden-Fowler equation (\ref{eq1}) with $p=3/2$. It is easy to check that in this case the following variables are homology invariant
\begin{equation}
    t=ax^{1/2}y^{1/6},\;\;\;u=by^{-4/3}y^\prime,
    \label{eq11}
\end{equation}
where $a$ and $b$ are some constants. We have
\begin{equation}
    x\frac{dt}{dx}=\frac{a}{2}x^{1/2}y^{1/6}+\frac{a}{6}x^{3/2}y^{-5/6}y^\prime=\frac{1}{2}t+\frac{a}{6}x^{3/2}y^{-5/6}y^\prime.
    \label{eq12}
\end{equation}
But
$$x^{3/2}y^{-5/6}y^\prime=x^{3/2}y^{3/6}\,y^{-8/6}y^\prime=\frac{t^3}{a^3}\,\frac{u}{b}.$$
Therefore,
\begin{equation}
    x\frac{dt}{dx}=\frac{t}{2}\left [1+\frac{1}{3a^2b}\,t^2u\right ].
    \label{eq13}
\end{equation}
On the other hand,
\begin{equation}
    x\frac{du}{dx}=-\frac{4}{3}bxy^{-7/3}y^{\prime\,2}+ bxy^{-4/3}y^{\prime\prime},
    \label{eq14}
\end{equation}
and
$$xy^{-7/3}y^{\prime\,2}=xy^{2/6}\,y^{-16/6}y^{\prime\,2}=\frac{t^2}{a^2}\,\frac{u^2}{b^2},\;\;\;
xy^{-4/3}y^{\prime\prime}=x^{1/2}y^{1/6}=\frac{t}{a}.$$
Therefore,
\begin{equation}
    x\frac{du}{dx}=\frac{bt}{a}\left [1-\frac{4}{3ab^2}\,tu^2\right ].
    \label{eq15}
\end{equation}
Dividing (\ref{eq15}) by (\ref{eq13}), we get the first order equation
\begin{equation}
    \frac{du}{dt}=2\frac{b}{a}\,\frac{1-\frac{4}{3ab^2}\,tu^2}{1+\frac{1}{3a^2b}\,t^2u}.
    \label{eq16}
\end{equation}
We can choose $a$ and $b$ to make this equation as simple as possible. Majorana's choice corresponds to
\begin{equation}
  \frac{4}{3ab^2}=1,\;\;\;\frac{1}{3a^2b}=-1,
  \label{eq17}
\end{equation}
which immediately gives
\begin{equation}
  a=12^{-1/3}=144^{-1/6},\;\;\;b=-4a=-\left(\frac{16}{3}\right )^{1/3}.
  \label{eq18}
\end{equation}
It is known that asymptotically, when $x\to \infty$, $y\to 144/x^3$. Then it is not difficult to see that Majorana's choice quite conveniently means that $t$ changes from $0$ to $1$ as $x$ changes from $0$ to $\infty$, and for $u$, as a function of $t$, we have $u(1)=1$. If you don't care about negative interval for $t$, an equally simple choice is $\frac{4}{3ab^2}=-1,\;\frac{1}{3a^2b}=1$. 

The Majorana transformation
\begin{equation}
    t=144^{-1/6}x^{1/2}y^{1/6},\;\;\;u=-\left(\frac{16}{3}\right )^{1/3}y^{-4/3}y^\prime,
    \label{eq19}
\end{equation}
which allows to reduce the Thomas-Fermi equation to the first order equation
\begin{equation}
    \frac{du}{dt}=-8\,\frac{1-tu^2}{1-t^2u},
    \label{eq20}
\end{equation}
seems mysterious when taken from nothing, like a magician pulling a hare out of a hat. However, for those initiated into the mystery of the trick, there is no magic here: the Majorana variables are just another set of homology-invariant variables.

\section{Dresner homology-invariant variables}
The final shade of the puzzle is why Majorana chose such homology-invariant variables. This choice does not seem trivial since the variables are hard to guess at first glance. We can give the following possible answer to this riddle.

In the case of the Thomas-Fermi model, the most simple and natural homology-invariant variables are \cite{Dresner}
\begin{equation}
    \tau=x^3y,\;\;\;s=x^4y^\prime.
    \label{eq21}
\end{equation}
Indeed, it is immediately evident that they are invariant under the homology transformations $x\to \lambda x$, $y\to \lambda^{-3}y$.
For these variables,
\begin{equation}
x\frac{d\tau}{dx}=3\tau+s,\;\;\;x\frac{ds}{dx}=4s+\tau^{3/2},
\label{eq21A}
\end{equation}
and therefore
\begin{equation}
    \frac{ds}{d\tau}=\frac{4s+\tau^{3/2}}{s+3\tau}.
    \label{eq22}
\end{equation}
We may want to get rid of irrationality of $\tau^{3/2}$ and introduce a new variables $t$ and $u$ via
\begin{equation}
     \tau=At^n,\;\;\;s=Bt^mu,
    \label{eq23}
\end{equation}
where $n$ is an even integer and $A$ and $B$ are some constants. Further note that $s^3/\tau^4=y^{\prime\,3}/y^4$ does not depend on the independent variable $x$. Therefore, it is natural to require that this ratio does not depend on the new independent variable $t$ either, and this is possible only when $4n-3m=0$.  The smallest integers $n$ and $m$ satisfying this relation for even $n$ are $n=6,\,m=8$. Then we get after some straightforward algebra
\begin{equation}
     \frac{ds}{d\tau}=\frac{B}{6A}\left [8t^2u+t^3\frac{du}{dt}\right ]=\frac{4Bt^2u+A^{3/2}t^3}{Bt^2u+3A}.
    \label{eq24}
\end{equation}
Solving this equation for $\frac{du}{dt}$, we get
\begin{equation}
    \frac{du}{dt}=\frac{6A}{B}\,\frac{A^{3/2}-\frac{4B^2}{3A}tu^2}{3A+Bt^2u}=2\,\frac{A^{3/2}}{B}\,\frac{1-\frac{4B^2}{3A^{5/2}}tu^2}{1+\frac{B}{3A}t^2u}.
    \label{eq25}
\end{equation}
It remains to choose the constants $A$ and $B$ in a convenient way, for which we require
\begin{equation}
    \frac{4B^2}{3A^{5/2}}=1,\;\;\;\frac{B}{3A}=-1.
    \label{eq26}
\end{equation}
The solution is $A=144,B=-3A$ and we recover the Majorana equation (\ref{eq20}) and the Majorana transformation
(\ref{eq19}). This is evident for $t=A^{-1/6}\tau^{1/6}=144^{-1/6}x^{1/2}y^{1/6}$, and for $u$ we can use
$$u=\frac{A^{4/3}}{B}\,\frac{s}{\tau^{4/3}}=-\frac{1}{3}A^{1/3}\,\frac{y^\prime}{y^{4/3}}=-\left(\frac{16}{3}\right )^{1/3}y^{-4/3}y^\prime.$$

To obtain the solutions $y(x)$, one must be able to transform back from the homology invariant variables to the original variables $(x,y)$. This can be done in the following way. Taking the logarithmic derivative with respect to $\tau$ of the first equation in (\ref{eq21}), we get
\begin{equation}
\frac{dy}{y}=\frac{d\tau}{\tau}-3\frac{dx}{x}.
\label{eq27}
\end{equation}
But {(\ref{eq21A})} indicates that $$\frac{dx}{x}=\frac{d\tau}{3\tau+s}.$$
Therefore
\begin{equation}
\frac{dy}{y}=\frac{sd\tau}{\tau(3\tau+s)}.
\label{eq28}
\end{equation}
In the original variables $(x,y)$, the one of boundary conditions is $y(0)=1$, and (\ref{eq21}) indicates that $\tau=0$, if $x=0$. Therefore, (\ref{eq28}) determines $y$ as a function of $\tau$, and then $x=(\tau/y)^{1/3}$ as a function of $\tau$, provided we know the solution $s(\tau)$ of the first order equation (\ref{eq22}):
\begin{equation}
y(\tau)=\exp{\left (\int\limits_0^\tau \frac{s(\sigma)\;d\sigma}{\sigma[3\sigma+s(\sigma)]}\right )},\;\;
x(\tau)=\tau^{1/3}\exp{\left (-\int\limits_0^\tau \frac{s(\sigma)\;d\sigma}{3\sigma[3\sigma+s(\sigma)]}\right )}.
\label{eq29}
\end{equation}
These relations constitute a parametric solution of the original Thomas-Fermi problem. We can easily express this parametric solution in terms of Majorana variables by using (\ref{eq23}):
\begin{equation}
y(t)=\exp{\left (\int\limits_0^t W(\sigma)d\sigma\right)}, \;\; x(t)=144^{1/3}t^2\exp{\left (-\frac{1}{3}\int\limits_0^t W(\sigma)d\sigma\right)},
\label{eq30}
\end{equation}
where
\begin{equation}
W(t)=-\frac{6tu(t)}{1-t^2u(t)}.
\label{eq31}
\end{equation}
Further details can be found in {\cite{Esposito2002}} and some subtletes are discussed in \cite{Seiler_2023}.

\section{Concluding remarks}
The method used by Majorana in 1928 to obtain a semi-analytical solution to the Thomas-Fermi equation remained unpublished and unknown until 2002 when it was revealed by Esposito in his illuminating pedagogical article  \cite{Esposito2002}.  The origin of Majorana's highly non-trivial change of variables may look enigmatic. In this note, we attempted to demystify the Majorana transformation using homology-invariant variables, a technique well known in the field of stellar astrophysics. In the light of this approach, the Majorana transformation becomes completely natural and transparent, and thanks to the simplicity of the method, one is tempted to assume that this is how Majorana found his transformation.

\section*{Aknowledgements}
We are grateful to the reviewers for useful comments.

\bibliography{MJT.bib}
\end{document}